\documentclass[a4paper,12pt]{article}

\usepackage{amsmath,amssymb,fourier}
\usepackage{graphicx}
\usepackage[french]{babel}

\newcommand{\cossi}{\,\includegraphics[viewport=0 0 9.7 8.8]{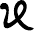}\,}
\newcommand{\cossii}{\,\includegraphics[viewport=0 2.5 4.6 10.5]{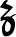}\,}
\newcommand{\cossiii}{\,\includegraphics[viewport=0 2.5 8.6 11.5]{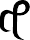}\,}
\newcommand{\cossiv}{\,\includegraphics[viewport=0 2.5 4.6 10.5]{coss2.pdf}\includegraphics[viewport=0 2.5 4.6 10.5]{coss2.pdf}\,}
\newcommand{\cossv}{\,\includegraphics[viewport=0 2 4.6 13]{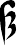}\,}
\newcommand{\cossvi}{\,\includegraphics[viewport=0 2.5 4.6 10.5]{coss2.pdf}\,\includegraphics[viewport=0 2.5 8.6 11.5]{coss3.pdf}\,}
\newcommand{\cossvii}{\,\includegraphics[viewport=0 2 15 13]{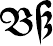}\,}
\newcommand{\cossviii}{\,\includegraphics[viewport=0 2.5 4.6 10.5]{coss2.pdf}\includegraphics[viewport=0 2.5 4.6 10.5]{coss2.pdf}\includegraphics[viewport=0 2.5 4.6 10.5]{coss2.pdf}\,}
\newcommand{\cossix}{\,\includegraphics[viewport=0 2.5 8.6 11.5]{coss3.pdf}\includegraphics[viewport=0 2.5 8.6 11.5]{coss3.pdf}\,}
\newcommand{\cossx}{\,\includegraphics[viewport=0 2.5 4.6 10.5]{coss2.pdf}\,\includegraphics[viewport=0 2 4.6 13]{coss5.pdf}\,}
\newcommand{\cossxi}{\mathfrak{C}\cossv\,}

\title{Quelques problèmes issus des \textit{Arithmetica philosophica} de Peter Roth}
\author{Erwan \textsc{Penchèvre} et David \textsc{Rabouin}}
\date{(\textit{draft})}

\begin{document}
\maketitle
Nous donnons ici traduction et commentaire de dix-sept problèmes issus du livre \textit{Arithmetica philosophica} de Peter Roth\footnote{Voir bibliographie en fin d'article, référence \cite{roth1608}.}. Tous ces problèmes font intervenir une équation algébrique à une inconnue, de degré 5 (\textit{surdesolid co\ss)}) ou 6 (\textit{zensicubic co\ss})~; en outre, quatorze d'entre eux font aussi intervenir une question de stéréométrie et une question concernant des «~nombres polygonaux~». Ces quatorze problèmes sont les seuls, dans les \textit{Arithmetica philosophica}, à faire intervenir ces trois thèmes à la fois au sein d'un même énoncé.

Chaque énoncé est précédé d'une référence. Par exemple «~XXIII f.~187~» désigne le problème~XXIII, folio~187. Les trois premiers problèmes nous permettent d'introduire certaines techniques. Le lecteur pressé pourra justement consulter d'emblée le problème~XXIII f.~187, car c'est un bon exemple.

\section{Exemples de binômes et résidus en \emph{surdesolid co\ss}}

\paragraph{I f. 184}
«~$\cossv-12\cossiv+10\cossiii=423\cossi-248\cossii+924$. \emph{Solution}~: $\cossi=6\pm\sqrt{3}$.~»

\paragraph{Commentaire}
Rappelons pour commencer le sens à donner aux symboles cossiques présents
dans les équations des \textit{Arithmetica philosophica}~:
\begin{align*}
  x=\cossi,\ x^2=\cossii,\ x^3=\cossiii,\ x^4=\cossiv,\ x^5=\cossv,\ x^6=\cossvi,\\
  x^7=\cossvii,\ x^8=\cossviii,\ x^9=\cossix,\ x^{10}=\cossx,\ x^{11}=\cossxi.
\end{align*}
Posons $P(x)=x^5-12x^4+10x^3+248x^2-423x-924$. Comme
il s'agit de trouver un binôme ou un résidu qui soit racine de ce
polynôme, on cherche un polynôme unitaire, facteur de degré 2 à
coefficients entiers de $P$.  Son terme constant doit diviser
$924=2^2.3.7.11$, on essaie tous les diviseurs possibles. Par exemple,
puisque $33$ est un diviseur, on cherche $y$ tel que $x^2-yx+33$
divise $P$. Pour ce faire, on pose la division euclidienne.  En
s'aidant d'une disposition en colonnes, on peut omettre de noter les
puissances de $x$~: les notations cossiques peuvent alors plutôt
servir à noter les puissances de l'indéterminée $y$. Pour faire que le
reste de la division s'annule identiquement, il faut trouver une
racine d'un polynôme en $y$~; comme on cherche une solution entière,
une telle racine est nécessairement un facteur du terme constant si le
polynôme en $y$ est unitaire. On trouvera $y=12$. Plus précisément
$P(x)=(x+4)(x^2-4x-7)(x^2-12x+33)$.  \emph{Autrement}~: on pourrait
d'abord avoir cherché un facteur linéaire évident (ici $x+4$) puis
avoir appliqué une des méthodes de résolution algébrique connues à
cette époque pour calculer les racines du polynôme quotient de degré~4
(méthode plus simple mais moins générale).

\paragraph{II f. 184}
«~$(\cossii+722)(\cossiii+3820)=2769560-50\cossiv-4416\cossi$. 
\emph{Solution}~: $\cossi=\sqrt{10}-2$.~»

\paragraph{Commentaire}
Voir commentaire précédent. L'équation est de la forme $P(x)=0$ avec
$P(x)=x^5+50x^4+722x^3+3820x^2+4416x-2^8.3^2.5$. On trouve facilement un
facteur quadratique car il y a trois facteurs linéaires, et 
$P(x)=(x+8)^2(x+30)(x^2+4x-6)$. L'unique racine positive du polynôme 
$x^2+4x-6$ est $\sqrt{10}-2$.

\paragraph{XXI f. 187}\label{questioXXI}
«~Soit un \textit{corpus irregulare} inscrit dans une sphère, à 18 faces
carrées et 8 faces triangulaires équilatérales régulières, tel que, quand il 
est posé sur un carré, la face supérieure est toujours aussi un carré. Quand
on le pose sur une base triangulaire, la face supérieure est toujours
parallèle à la face inférieure (c'est-à-dire la base) et c'est aussi un
triangle. L'arête de ce corps à vingt-quatre sommets mesure tant d'unités que
la somme des deux vraies racines cinquièmes de 
$71\cossiii+75\cossi-2\cossiv-78\cossii-450$.
Combien le diamètre de la sphère circonscrite mesure-t-il~? 
\emph{Solution}~: $\sqrt{405+\sqrt{52488}}$.~»

\paragraph{Commentaire}
Ce polyèdre semi-régulier obtenu par double troncature à partir d'un 
cube ou d'un octaèdre est décrit par Stévin\footnote{Voir \cite{stevin1583}, 
  définition 15 p.~51.} qui le nomme \textit{bistruncatum cubum primum}. Il en
dessine le patron et calcule le rapport entre son arête et le rayon de la sphère 
circonscrite\footnote{Voir \cite{stevin1583}, p.~75 et 58.}.
Ce polyèdre est aussi étudié par Descartes 
dans \cite{descartes}, comme tous les polèdres semi-réguliers obtenus par
troncature~: dans le commentaire de Costabel, il est noté
$18[4]\ 8[3]$. Soit $d$ la longueur d'une arête et $\rho$ le rayon de 
la sphère circonscrite, alors $(2\rho)^2=(5+2\sqrt{2})d^2$. On a $d=x_0+x_1$ 
où $x_0$ et $x_1$ sont deux solutions réelles positives de l'équation
$$x^5=71x^3+75x-2x^4-78x^2-450.$$
Posons $P(x)=x^5+2x^4-71x^3+78x^2-75x+450$. On a $450=2.3^2.5^2$, on cherche
un facteur quadratique de la forme $x^2-yx+15$. On trouve que
$$P(x)=(x+10)(x^2+x+3)(x^2-9x+15).$$
Ses deux seules racines positives sont $\dfrac{9\pm\sqrt{21}}{2}$, donc
$d=9$ et $2\rho=9\sqrt{5+2\sqrt{2}}=\sqrt{405+\sqrt{52488}}$. \emph{Remarque}~:
si l'on a une idée préalable de l'ordre de grandeur de la solution, si
l'on sait par exemple que $d<10$, comme la somme des deux racines du
facteur quadratique doit être inférieure à dix, leur produit ne pourra être
supérieur à 25. Ainsi, le terme constant du facteur quadratique devra être
un diviseur de 450 \emph{inférieur à 25}. On a le choix entre
1, 2, 3, 5, 6, 9, 10, 15 ou 18. Dans tous les problèmes de ce genre, s'il
s'agit seulement de \emph{retrouver} les calculs menant à une solution déjà
connue (puisque Roth donne lui-même la solution), on peux considérablement
abréger la recherche en utilisant cette remarque.

\paragraph{XXII f. 187}
«~Soit une sphère dont le diamètre mesure autant d'unités que la plus grande 
vraie racine \textit{pentachilitetracosihexagonale} de 
$\cossv+3550\cossii-1817\cossi+5232-4\cossiv-84\cossiii$.
On y inscrit un \textit{corpus} irrégulier à 6 faces carrées et 8 faces
hexagonales équilatérales régulières, tel que, quand il est posé sur un carré,
la face supérieure est toujours aussi un carré. Quand on le pose sur
une base hexagonale, la face supérieure est toujours parallèle à la face
inférieure (c'est-à-dire la base) et c'est aussi un hexagone. Combien
une arête de ce corps à vingt-quatre sommets inscrit dans la sphère 
donnée mesure-t-elle~? \emph{Solution}~: 
$\sqrt{3+\dfrac{3}{5}}+\sqrt{\dfrac{1}{5}}$.~»
 
\paragraph{Commentaire} L'équation donnée dans cet énoncé est probablement
erronée. C'est fâcheux, car il s'agit du premier d'une série de problèmes
faisant intervenir à la fois des polyèdres, des équations algébriques, et 
des nombres figurés plans (le lecteur comprendra plus aisément s'il
consulte d'abord le problème~XXIII). Le polyèdre $8[6]\ 6[4]$ ici en question
est obtenu par troncature d'un octaèdre, les arêtes étant coupées au tiers. 
Stévin dessine le patron de cet \textit{octoedrum truncatum per laterum
tertias}\footnote{Voir \cite{stevin1583} p.~78.}, et il calcule $(2\rho)^2=10d^2$~; 
d'après la solution donnée par Roth, il faut donc que 
$$x=\sqrt{10}d=\sqrt{10}\left(\sqrt{3+\frac{3}{5}}+\sqrt{\frac{1}{5}}
\right)=6+\sqrt{2}.$$
Le polynôme quadratique unitaire dont ce nombre est la plus grande racine 
positive est $x^2-12x+34$. Pour comprendre l'équation décrite dans l'énoncé,
il faut savoir qu'un nombre 5406-gonal est un nombre de la forme
$$O_n(x)=\frac{n-2}{2}x^2-\frac{n-4}{2}x,$$
avec $n=5406$~; l'équation de l'énoncé s'écrit donc
$$2702x^2-2701x=x^5+3550x^2-1817x+5232-4x^4-84x^3.$$
Elle est de la forme $P(x)=0$ avec
$$P(x)=x^5-4x^4-84x^3+848x^2+884x+5232.$$
Ce polynôme n'est certes pas divisible par $x^2-12x+34$~; il y a donc 
certainement une coquille dans l'équation (on ne veut pas mettre en cause la 
solution donnée par Roth ni les autres données de l'énoncé, car on retrouve 
le même facteur quadratique dans les problèmes suivants). \emph{Remarque}~: on
s'abstiendra en général de chercher à corriger de telles erreurs d'énoncé
sauf les plus évidentes~; mais changer \emph{deux} coefficients suffit toujours
à rendre un polynôme quelconque divisible par un facteur quadratique donné.
Ici, on pourrait proposer le polynôme $P$ suivant en remplacement de celui de
l'énoncé~:
$$P(x)=x^5-84x^3+96x^2+884x=x(x^2+12x+26)(x^2-12x+34).$$

\paragraph{XXIII f. 187}
«~Soit une autre sphère dont le diamètre est égal à la plus petite vraie racine
\textit{triacosiaicosihexagonale} de 
$\cossv+578\cossii+7072-2861\cossi-20\cossiii-8\cossiv$. On
y inscrit un \textit{corpus} irrégulier à 6 faces carrées et 8 faces
triangulaires équilatérales régulières, tel que, quand il est posé sur un
carré, la face supérieure est toujours aussi un carré. Quand on le pose sur
une base triangulaire, la face supérieure lui est toujours parallèle et 
c'est aussi un triangle. Combien une arête de ce corps à douze sommets inscrit
dans la sphère donnée mesure-t-elle~? \emph{Solution}~: 
$3-\sqrt{\dfrac{1}{2}}$.~»

\paragraph{Commentaire}
Le polyèdre $6[4]\ 8[3]$ est obtenu par troncature à partir d'un cube ou
d'un octaèdre, les arêtes étant coupées à la moitié. Stévin dessine le patron
de cet \textit{octoedrum truncatum per laterum media}\footnote{C'est ainsi que Stévin le nomme dans \cite{stevin1583} p.~77. Cf. \textit{cuboctahedron}, \cite{coxeter1948} \S~2.3 p.~18.}, et il montre que $\rho=d$. Or l'énoncé demande que $\rho$ soit la 
plus petite solution positive de l'équation
$$O_{326}(x)=x^5+578x^2+7072-2861x-20x^3-8x^4,$$
où $O_{326}(x)$ est le nombre 326-gonal de côté $x$ défini par
$$O_n(x)=\frac{n-2}{2}x^2-\frac{n-4}{2}x,$$
avec $n=326$. L'équation est de la forme $P(x)=0$ avec
$$P(x)=x^5-8x^4-20x^3+416x^2-2700x+7072.$$
On a $7072=2^5.13.17$, on trouve un facteur quadratique de la forme
$x^2-yx+34$ et
$$P(x)=(x+8)(x^2-4x+26)(x^2-12x+34).$$
Ses seules racines positives sont $6\pm\sqrt{2}$, donc $2\rho=6-\sqrt{2}$ et
$d=3-\sqrt{\dfrac{1}{2}}$. 

\paragraph{XXIV f. 187}\label{questioXXIV}
«~Soit un \textit{corpus irregulare} (inscrit dans une sphère creuse) à
20 faces hexagonales et douze faces pentagonales régulières, tel que, quand
il est posé sur un pentagone, la face supérieure est toujours aussi un
pentagone. Quand on le pose sur une base hexagonale, la face supérieure
lui est toujours parallèle et c'est aussi un hexagone. La longueur de 
l'arête de ce corps à soixante sommets est la somme des deux vraies racines
\textit{heptacosihexacontapentagonales} de $\cossv
+\left(689+\frac{1}{2}\right)\cossii+6188-9\cossiv-4\cossiii
-\left(2632+\frac{1}{2}\right)\cossi$.
Combien le diamètre de la sphère circonscrite mesure-t-il~? 
\emph{Solution}~: le diamètre mesure $\sqrt{2088+\sqrt{2099520}}$.

«~Le grand \textit{Mathematicus} Simon Jacob, \textit{Rechenmeister} à
Francfort offre un exemple de ce \textit{corpus} dans son \textit{
Arithmetica} (\textit{in quarto, folio 338}). J'ai résolu son problème
de bien des manières~; un jour si Dieu le veut, je publierai cela avec
d'autres questions sur les polyèdres issues de son livre.~»

\paragraph{Commentaire}
Le polyèdre $20[6]\ 12[5]$ est obtenu par troncature à partir de l'icosaèdre,
les arêtes étant coupées au tiers. Stévin dessine le patron de ce
\textit{truncatum icosaedrum per laterum tertias}\footnote{Voir \cite{stevin1583} p.~81.}
qu'il affirme ne pas avoir rencontré dans les ouvrages de ses 
prédécesseurs\footnote{Voir \cite{stevin1583} p.~47.}, et il calcule 
$(2\rho)^2=\dfrac{29+9\sqrt{5}}{2}d^2$.
L'énoncé demande que $d=x_0+x_1$ où $x_0$ et $x_1$ sont deux racines réelles
positives de l'équation
$$O_{765}(x)=x^5+\left(689+\frac{1}{2}\right)x^2+6188-9x^4-4x^3
-\left(2632+\frac{1}{2}\right)x.$$
Cette équation est de la forme $P(x)=0$ avec
$$P(x)=x^5-9x^4-4x^3+308x^2-2252x+6188.$$
On a $6188=2^2.7.13.17$, et on trouve un facteur quadratique de la forme
$x^2-yx+34$~:
$$P(x)=(x+7)(x^2-4x+26)(x^2-12x+34).$$
Les seules racines positives de $P$ sont $6\pm\sqrt{2}$, donc $d=12$ et
$$2\rho=12\sqrt{\frac{29+9\sqrt{5}}{2}}=\sqrt{2088+\sqrt{2099520}}.$$

\paragraph{XXV f. 187}
«~Soit un \textit{corpus irregulare} à 12 faces pentagonales et 20 faces
triangulaires régulières équilatérales, tel que, quand on le pose sur
un pentagone, cette base est parallèle à la face supérieure qui est 
toujours aussi un pentagone. Quand on le pose sur un triangle, alors
la face supérieure est toujours aussi un triangle. Le diamètre d'un cercle
circonscrivant l'un des pentagones de ce \textit{corpus} à trente sommets
est égal à la différence entre les deux vraies racines sixièmes de
$1\cossvi+1\cossi+1172\cossii+13260-1\cossiv-132\cossiii-5836\cossi$.
Combien le diamètre de la sphère circonscrite mesure-t-il~?
\emph{Solution}~: $\sqrt{20+\sqrt{80}}$.~»

\paragraph{Commentaire}
Le polyèdre $20[3]\ 12[5]$ est obtenu par troncature à partir de l'icosaèdre
ou du dodécaèdre, les arêtes étant coupées à la moitié. Stévin dessine le
patron de ce \textit{truncatum icosaedrum per laterum media}\footnote{Voir \cite{stevin1583} p.~80. Cf. \textit{icosidodecahedron}, \cite{coxeter1948} \S~2.3 p.~18.} qu'il affirme ne pas avoir rencontré dans les
ouvrages de ses prédécesseurs, et il calcule $(2\rho)^2=(1+\sqrt{5})^2d^2$.
Notons $r$ le rayon du cercle circonscrit à un pentagone de côté $d$,
on a
$$d=r\sqrt{\frac{5-\sqrt{5}}{2}}.$$
L'énoncé demande que $2r=x_0-x_1$ où $x_0$ et $x_1$ sont deux racines
réelles positives de l'équation $P(x)=0$ avec
$$P(x)=x^5-x^4-132x^3+1172x^2-5836x+13260.$$
On a $13260=2^2.3.5.13.17$, et on trouve un facteur quadratique (encore le
même) de la forme $x^2-yx+34$~:
$$P(x)=(x+15)(x^2-4x+26)(x^2-12x+34).$$
Les seules racines positives de $P$ sont $6\pm\sqrt{2}$, donc $2r=2\sqrt{2}$
et $d=\sqrt{5-\sqrt{5}}$. Finalement
$$2\rho=(1+\sqrt{5})\sqrt{5-\sqrt{5}}=\sqrt{20+\sqrt{80}}.$$

\paragraph{XXVI f. 188}
«~Soit une sphère dont le diamètre est égal à la plus grande vraie racine
\textit{enneacositriacontahexagonale} de $4\cossiv+\left(29+\frac{1}{2}
\right)\cossiii-\left(131+\frac{5}{16}\right)\cossi-1\cossv
+\left(371+\frac{1}{2}\right)\cossii-\left(1141+\frac{7}{8}\right)$
On inscrit dans cette sphère un \textit{corpus} irrégulier à quatre faces
hexagonales et quatre faces triangulaires équilatérales, tel que, quand on
le pose sur un hexagone, la face supérieure est toujours parallèle à la base
et elle est triangulaire. Inversement, quand on le pose sur un triangle,
la face supérieure est toujours un hexagone. Combien l'arête de ce
\textit{corpus} à douze sommets inscrit dans la sphère mesure-t-elle~?
\emph{Solution}~: $\sqrt{1+\dfrac{3}{22}}+\sqrt{\dfrac{13}{88}}$. Il ne
peut en être autrement.~»

\paragraph{Commentaire}
Le polyèdre $4[6]\ 4[3]$ est obtenu par troncature à partir du tétraèdre,
les arêtes étant coupées au tiers.  Stévin dessine le patron de ce 
\textit{truncatum tetraedrum}\footnote{Voir \cite{stevin1583} p.~73.}, et il calcule 
$(2\rho)^2=\dfrac{11}{2}d^2$. L'énoncé demande que $2\rho=x$ où
$x$ est la plus grande racine réelle positive de l'équation
$$O_{936}(x)=4x^4+\left(29+\frac{1}{2}\right)x^3
-\left(131+\frac{5}{16}\right)x-x^5+\left(371+\frac{1}{2}\right)x^2
-\left(1141+\frac{7}{8}\right).$$
Cette équation est de la forme $P(x)=0$ avec
$$P(x)=x^5-4x^4-\left(29+\frac{1}{2}\right)x^3+\left(95+\frac{1}{2}\right)x^2
-\left(334+\frac{11}{16}\right)x+\left(1141+\frac{7}{8}\right).$$
Les \textit{Rechenmeister} commençaient alors habituellement par faire
un changement de variable pour transformer le polynôme en un polynôme
unitaire à coefficients entiers~; posons donc $x=\dfrac{z}{2}$. On a
$$32P\left(\frac{z}{2}\right)=z^5-8z^4-118z^3+764z^2-5355z+36540.$$
On a $36540=2^2.3^2.5.7.29$, et on trouve un facteur quadratique de la
forme $z^2-yz+87$~:
$$32P\left(\frac{z}{2}\right)=(z+12)(x^2+35)(z^2-20z+87).$$
La plus grande racine réelle positive de ce polynôme est $z=10+\sqrt{13}$,
donc $2\rho=x=5+\dfrac{\sqrt{13}}{2}$, et
$$d=\left(5+\frac{\sqrt{13}}{2}\right)\sqrt{\frac{2}{11}}=2\times\left(
\sqrt{1+\frac{3}{22}}+\sqrt{\frac{13}{88}}\right).$$
La solution donnée par Roth est donc erronée d'un facteur $\frac{1}{2}$,
certainement dû à une confusion entre rayon et diamètre de la 
sphère\footnote{\textbf{Il faudra vérifier que l'erreur ne vient pas de Stévin…}} 
(ce ne peut être dû au changement de variable 
car une telle erreur aurait conduit à un facteur 2 et non $\frac{1}{2}$).

\paragraph{XXVII f. 188}
«~Soit un \textit{corpus irregulare} à six faces carrées équilatérales et
régulières, et douze faces triangulaires non équilatérales. On suppose que
chaque hauteur abaissée d'un sommet d'une telle face triangulaire sur 
l'arête d'une face carrée (en son milieu) est égale à l'arête des faces
carrées de ce \textit{corpus}. Quand on le pose sur une face carrée, la
face supérieure est toujours aussi un carré. Quand on le pose sur une face
triangulaire, la face supérieure est toujours parallèle à la base, et
c'est aussi un triangle. Le diamètre de la sphère circonscrite (telle que ce
\textit{corpus} la touche en chaque sommet) est égal à la plus petite vraie
racine \textit{chiliadiacosiagonale} de $1\cossv+\left(755+\frac{1}{2}
\right)\cossii+\left(1522+\frac{1}{2}\right)-2\cossiv-\left(49+\frac{1}{2}
\right)\cossiii-\left(1207+\frac{11}{16}\right)\cossi$. Combien l'arête 
d'une face carrée de ce \textit{corpus} à quatorze sommets mesure-t-elle~? Et 
les arêtes des triangles~? \emph{Solution}~: les arêtes des carrés font
$\sqrt{5}-\sqrt{\dfrac{13}{20}}$, et les deux autres arêtes de chaque triangle
font $2+\dfrac{1}{2}-\sqrt{\dfrac{13}{16}}$.~»

\paragraph{Commentaire}
Le polyèdre décrit ici appartient à la famille des «~diamants hexagonaux
allongés~». Ce n'est pas un polyèdre semi-régulier. Roth commet ici une erreur 
en affirmant qu'il existe une sphère circonscrite «~telle qu'il la touche en 
chaque sommet~» (\textit{{da\ss } es mit allen seinen spitzen darin anstöst}).
Il existe en effet une unique sphère contenant les douze sommets des 
faces carrées, mais les deux sommets restant n'y sont pas. Notons $\rho$ le
rayon de cette sphère, $d_1$ l'arête des faces carrées, et $d_2$ la longueur
de chacune des deux autres arêtes des faces triangulaires. \`A cause de la
donnée concernant les hauteurs des faces triangulaires, on a
$$d_2=d_1\frac{\sqrt{5}}{2}.$$
On montre d'autre part que $(2\rho)^2=5d_1^2$. Or l'énoncé affirme que
$2\rho$ est la plus petite racine positive de l'équation\footnote{
Exceptionnellement, on s'est ici permis de corriger une coquille dans un des
coefficient de l'équation~: au lieu de 1107, il est en effet écrit 1207. Dans
l'exemplaire que nous avons consulté (British Library), un lecteur avait
d'ailleurs apporté cette correction en marge.}
$$O_{1200}(x)=x^5+\left(755+\frac{1}{2}\right)x^2
+\left(1522+\frac{1}{2}\right)-2x^4-\left(49+\frac{1}{2}\right)x^3
-\left(1107+\frac{11}{16}\right)x.$$
Cette équation est de la 
forme $P(x)=0$ avec
$$P(x)=x^5-2x^4-\left(49+\frac{1}{2}\right)x^3
+\left(156+\frac{1}{2}\right)x^2-\left(509+\frac{11}{16}\right)x
+\left(1522+\frac{1}{2}\right).$$
On pose $x=\dfrac{z}{2}$~; alors
$$32P\left(\frac{z}{2}\right)=z^5-4z^4-198z^3+1252z^2-8155z+48720.$$
On a $48720=2^4.3.5.7.29$, et on trouve un facteur quadratique de la
forme $z^2-yz+87$ (le même que dans le problème précédent)~:
$$32P\left(\frac{z}{2}\right)=(z+16)(z^2+35)(z^2-20z+87).$$
La plus petite racine réelle positive est $2\rho=5-\dfrac{\sqrt{13}}{2}$~;
on trouve donc bien 
$$d_1=\sqrt{5}-\sqrt{\dfrac{13}{20}},\text{ et }d_2=2+\frac{1}{2}
-\sqrt{\frac{13}{16}}.$$

\paragraph{XXVIII f. 188}
«~Soit un \textit{corpus irregulare} à 6 faces octogonales et 8 faces
triangulaires équilatérales régulières, tel que, quand on le pose sur un
triangle, la face supérieure est toujours aussi un triangle. Quand on le pose
sur un octogone, la face supérieure est toujours parallèle à la face
inférieure et elle est aussi un octogone. Le diamètre de la sphère
circonscrite à ce \textit{corpus} est la somme des deux vraies racines
\textit{hebdomicontahenagonales} de $1\cossv+1\cossiv+\left(282
+\frac{1}{2}\right)\cossii+\left(2093+\frac{7}{16}\right)
-\left(79+\frac{1}{2}\right)\cossiii-\left(805+\frac{11}{16}\right)\cossi$.
Combien l'arête de ce \textit{corpus} à vingt-quatre sommets mesure-t-elle~?
\emph{Solution}~: $\sqrt{\dfrac{700-\sqrt{320000}}{17}}$.~»

\paragraph{Commentaire}
Le polyèdre $6[8]\ 8[3]$ est obtenu par troncature à partir du cube, les
arêtes étant coupées au tiers. Stévin dessine le patron de ce 
\textit{truncatum cubum per laterum divisionis in tres partes}\footnote{Voir
\cite{stevin1583} p.~74.} et il calcule $(2\rho)^2=(7+4\sqrt{2})d^2$.
L'énoncé demande que $2\rho=x_0+x_1$ où $x_0$ et $x_1$ sont deux racines
réelles positives de l'équation
$$O_{71}(x)=x^5+x^4+\left(282+\frac{1}{2}\right)x^2
+\left(2093+\frac{7}{16}\right)-\left(79+\frac{1}{2}\right)x^3
-\left(805+\frac{11}{16}\right)x.$$
Cette équation est de la forme $P(x)=0$ avec
$$P(x)=x^5+x^4-\left(79+\frac{1}{2}\right)x^3+248x^2-\left(772+\frac{3}{16}
\right)x+\left(2093+\frac{7}{16}\right).$$
On pose $x=\dfrac{z}{2}$ et on a
$$32P\left(\frac{z}{2}\right)=z^5+2z^4-318z^3+1984z^2-12355z+66990.$$
On a $66990=2.3.5.7.11.29$ et on trouve encore un facteur quadratique de la 
forme $z^2-yz+87$~:
$$32P\left(\frac{z}{2}\right)=(z+22)(z^2+35)(z^2-20z+87).$$
Finalement les deux seules racines réelles positives de $P$ sont
$5\pm\dfrac{\sqrt{13}}{2}$, et $2\rho=10$, donc on a bien
$$d=\frac{2\rho}{\sqrt{7+4\sqrt{2}}}=\sqrt{\frac{700-400\sqrt{2}}{17}}.$$

\paragraph{XXIX f. 188}
«~Soit un autre \textit{corpus irregulare}, inscrit dans une sphère, à
20 faces triangulaires et 12 faces décagonales équilatérales régulières,
tel que, quand on le pose sur un triangle, la face supérieure est toujours
aussi un triangle. Quand on le pose sur un décagone, la face supérieure est
toujours parallèle à la face inférieure est c'est aussi un décagone.
L'arête de ce \textit{corpus} à soixante sommets fait la différence entre les
deux vraies racines \textit{tessaracontaoctogonales} de $1\cossv+3\cossiv
+332\cossii+\left(2474+\frac{1}{16}\right)-\left(99+\frac{1}{2}\right)
\cossiii-\left(969+\frac{3}{16}\right)\cossi$. Combien le diamètre de la
sphère circonscrite mesure-t-il~? \emph{Solution}~: $\sqrt{1850
+\sqrt{2812500}}$.~»

\paragraph{Commentaire}
Le polyèdre $12[10]\ 20[3]$ est obtenu par troncature à partir du dodécaèdre,
les arêtes étant coupées au tiers. Stévin dessine le patron de ce
\textit{truncatum dodecaedrum per laterum divisiones in tres partes}\footnote{Voir \cite{stevin1583} p.~79.} qu'il affirme ne pas avoir rencontré dans les
ouvrages de ses prédécesseurs, et il calcule
$(2\rho)^2=\dfrac{37+15\sqrt{5}}{2}d^2$.
L'énoncé demande de trouver les racines positives de l'équation
$$O_{48}(x)=x^5+3x^4+332x^2+\left(2474+\frac{1}{16}\right)
-\left(99+\frac{1}{2}\right)x^3-\left(969+\frac{3}{16}\right)x.$$
Cette équation est de la forme $P(x)=0$ avec
$$P(x)=x^5+3x^4-\left(99+\frac{1}{2}\right)x^3+309x^2-\left(947
+\frac{3}{16}\right)x+\left(2474+\frac{1}{16}\right).$$
On pose $x=\dfrac{z}{2}$~:
$$32P\left(\frac{z}{2}\right)=z^5+6z^4-398z^3+2472z^2-15155z+79170.$$
On a $79170=2.3.5.7.13.29$ et on trouve encore le même facteur quadratique,
de la forme $z^2-yz+87$~:
$$32P\left(\frac{z}{2}\right)=(z+26)(z^2+35)(z^2-20z+87).$$
Les deux seules racines positives de $P$ sont donc 
$x_0=5-\dfrac{\sqrt{13}}{2}$ et $x_1=5+\dfrac{\sqrt{13}}{2}$. Si, comme
le demande l'énoncé, $d=x_1-x_0$, on obtient $d=\sqrt{13}$ qui conduit à
un résultat nettement inférieur à la solution donnée par Roth~; mais la valeur
donnée par Roth s'explique aisément, car poser $d=x_0+x_1=10$ conduit 
précisément à cette solution~:
$$2\rho=10\sqrt{\frac{37+15\sqrt{5}}{2}}=\sqrt{1850+\sqrt{2812500}}.$$
Ici, il faut donc penser à une erreur de signe dans la lecture de l'énoncé 
de la part de celui qui a calculé la solution donnée par Roth.

\paragraph{XXX f. 188}
«~Soit un \textit{corpus irregulare} à six faces octogonales, huit faces
hexagonales et douze faces carrées équilatérales régulières, tel que,
posé sur une octogone, la face supérieure est aussi un octogone~; posé
sur un carré, la face supérieure est toujours aussi un carré~; posé sur
un hexagone, la face supérieure est aussi un hexagone~; et les faces
supérieures sont toujours parallèles à la base. L'arête de ce \textit{corpus}
est la somme des deux vraies racines \textit{diacositriacontadygonales}
de $1\cossv+137\cossii+1848-1\cossiv-42\cossiii-382\cossi$. Combien
le diamètre de la sphère circonscrite à ce \textit{corpus} à quarante-huit
sommets mesure-t-il~? \emph{Solution}~:
$$\sqrt{\frac{\sqrt{1400000000+\sqrt{320000000000000000}}+10000}
{\sqrt{35000+\sqrt{200000000}}-125-\sqrt{5000}}}$$
ou bien
$$\sqrt{\sqrt{2240000+\sqrt{819200000000}}+400}\text{.~»}$$

\paragraph{Commentaire}
Le polyèdre $12[4]\ 8[6]\ 6[8]$ est obtenu par double troncature à partir
du cube~; Descartes le connaît mais n'a pas calculé le rapport entre son 
arête et le diamètre de la sphère circonscrite. Stévin dessine le patron
de ce polyèdre qu'il appelle \textit{bistruncatum cubum secundum}\footnote{Voir \cite{stevin1583} p.~76.}, et il calcule $(2\rho)^2=(13+6\sqrt{2})d^2$ 
(\textit{distinctio} 14 p.~58, et figure p.~61). L'énoncé demande que
$d=x_0+x_1$ où $x_0$ et $x_1$ sont deux racines réelles positives de
l'équation
$$O_{232}(x)=x^5+137x^2+1848-x^4-42x^3-382x.$$
Cette équation est de la forme $P(x)=0$ avec
$$P(x)=x^5-x^4-42x^3+22x^2-268x+1848.$$
On a $1848=2^3.3.7.11$ et on trouve un facteur quadratique de la forme
$x^2-yx+22$~:
$$P(x)=(x+7)(x^2+2x+12)(x^2-10x+22).$$
Les seules racines réelles positives sont $5\pm\sqrt{3}$, donc $d=10$ et~:
$$2\rho=10\sqrt{13+6\sqrt{2}}\simeq 46{,}35.$$

Reste à expliquer les deux expressions compliquées données par Roth comme
solutions. Elles sont erronées. Après simplification des radicaux et des 
fractions, on trouve que ces deux expressions sont respectivement égales à~:
$$20\sqrt{\frac{1+\sqrt{14+4\sqrt{2}}}{2\sqrt{14+4\sqrt{2}}-5-2\sqrt{2}}}$$
et
$$20\sqrt{1+\sqrt{14+4\sqrt{2}}}.$$
Leurs valeurs arrondies au centième près sont respectivement $45,74$ et
$46,62$. Ces deux valeurs étant distinctes, on peut penser à une erreur, non
pas dans l'énoncé ou l'équation proposée\footnote{L'équation est hors de
doute car le problème suivant contient le même facteur quadratique…},
mais bien dans la solution (probablement dans la construction de la sphère 
circonscrite, problème délicat auquel Stévin avait d'ailleurs dédié une 
figure séparée).

\section{Exemple et questions en \textit{zensicubic co\ss}}

\paragraph{I f. 189}\label{questioI}
«~Soit un autre \textit{corpus} irrégulier, à vingt faces triangulaires,
trente faces carrées et douze faces pentagonales, tel que, quand on le pose
sur un triangle, la face supérieure est toujours parallèle à la base sur
laquelle il est posé, et c'est aussi un triangle. Quand on le pose sur un
carré, il y a la même propriété~: la face supérieure est aussi un carré.
De même pour les faces pentagonales~: quand on le pose sur l'une d'elles, la
face supérieure est toujours aussi un pentagone. Le diamètre de la sphère
circonscrite à ce \textit{corpus} à soixante sommets fait la différence des
deux vraies racines sixièmes de $54\cossiv+482\cossiii+4\cossii+1368\cossi
-11\cossv-22176$. Combien l'arête de ce \textit{corpus} mesure-t-elle~?
\emph{Solution}~: $\sqrt{\dfrac{132-\sqrt{11520}}{41}}$.~»

\paragraph{Commentaire}
Le polyèdre $20[3]\ 12[5]\ 30[4]$ est obtenu par double
troncature à partir de l'icosaèdre ou du dodécaèdre. 
Stévin ne le mentionne pas dans
\cite{stevin1583}. On l'appelle aujourd'hui «~rhombicosidodécaèdre~».
On peut montrer que\footnote{Les coordonnées des sommets d'un rhombicosidodécaèdre
  d'arête $d=2\varphi^{-1}$, où $\varphi=\dfrac{1+\sqrt{5}}{2}$, sont données par Coxeter \cite{coxeter1948} tabl.~V (iv) p.~299. Ces coordonnées permettent de calculer aisément \(\left(\dfrac{2\rho}{d}\right)^2=5\varphi^2+\varphi^4=7+8\varphi=11+4\sqrt{5}\).} $(2\rho)^2=(11+4\sqrt{5})d^2$. 
L'énoncé demande que $2\rho=x_0-x_1$ où $x_0>x_1$ sont
deux racines positives de l'équation $P(x)=0$, avec
$$P(x)=x^6+11x^5-54x^4-482x^3-4x^2-1368x+22176.$$
On a $22176=2^5.3^2.7.11$, et on trouve un facteur quadratique
de la forme $x^2-yx+22$ (le même que dans le problème précédent)~:
$$P(x)=(x+7)(x+12)(x^2+2x+12)(x^2-10x+22).$$
Les deux seules racines réelles positives sont $5\pm\sqrt{3}$, donc
$2\rho=2\sqrt{3}$, et $d=\sqrt{\dfrac{132-\sqrt{11520}}{41}}$ comme annoncé.

\paragraph{II f. 189}
«~Il existe encore un ingénieux et surprenant \textit{corpus} irrégulier à six 
faces carrées et 32 faces triangulaires équilatérales régulières, tel que,
quand on le pose sur un carré, la face supérieure lui est toujours parallèle
et c'est aussi un carré (mais ils sont de biais l'un par rapport à l'autre, de
sorte qu'une perpendiculaire abaissée d'un des sommets de la face carrée
supérieure sur la face carrée inférieure ne tombera pas loin du milieu d'une
arête de cette face inférieure). Quand on le pose sur une face triangulaire 
à laquelle sont attachées trois faces triangulaires semblables, la face
supérieure est toujours parallèle à la face inférieure et elle est aussi
triangulaire~; lui sont aussi attachées trois faces triangulaires. Chaque
arête de ce \textit{corpus} à vingt-quatre sommets est égale à la somme des
deux vraies racines \textit{dischilipentacositriacontahexagonales} de
$1\cossvi+11\cossv+\left(2911+\frac{5}{16}\right)\cossii
+\left(9756+\frac{1}{2}\right)-\left(81+\frac{1}{2}\right)\cossiv
-573\cossiii-\left(4235+\frac{15}{16}\right)\cossi$. Combien le diamètre
de la sphère circonscrite mesure-t-il~? \emph{Solution}~: le diamètre
de la sphère mesure
$$\scriptstyle\sqrt{\frac{2666+\frac{2}{3}-\sqrt[3]{703648148+\frac{4}{27}
-\sqrt{407329221193415637+\frac{209}{243}}}-\sqrt[3]{703648148+\frac{4}{27}
+\sqrt{407329221193415637+\frac{209}{243}}}}{16+\frac{2}{3}-\sqrt[3]{703
+\frac{35}{54}-\sqrt{407329+\frac{215}{972}}}-\sqrt[3]{703+\frac{35}{54}
+\sqrt{407329+\frac{215}{972}}}}}\displaystyle\text{.~»}$$

\paragraph{Commentaire}
Le polyèdre $6[4]\ 32[3]$ ne peut être obtenu par troncature à partir des
polyèdres réguliers. Stévin en détient la connaissance de Dürer, et il en 
dessine le patron\footnote{Voir \cite{stevin1583} p.~83.}. Kepler le nommera \textit{cubus
simus}. Puisque c'est le dernier polyèdre semi-régulier décrit par Roth
dans les \textit{Arithmetica philosophica}, il faut remarquer qu'il en
omet donc deux~:
\begin{itemize}
\item un polyèdre $30[4]\ 20[6]\ 12[10]$, obtenu par double troncature à 
partir de l'icosaèdre ou du dodécaèdre~; Descartes le connaît mais n'a pas 
calculé le rapport entre son arête et le diamètre de la sphère circonscrite~;
Stévin ne le mentionne pas dans \cite{stevin1583}~;
\item un polyèdre $12[5]\ 80[3]$, qui ne peut être obtenu par troncature 
à partir des polyèdres réguliers, mais qui est au dodécaèdre ce que le 
\textit{cubus simus} est au cube~; Descartes ne mentionne pas ce polyèdre, 
Stévin non plus (\textbf{et donc Dürer non plus, je suppose}).
\end{itemize}
Revenons au \textit{cubus simus}. Stévin écrivait en 1583~: «~Faute
de savoir comment l'obtenir par troncature, c'est-à-dire d'en 
  connaître la véritable origine, nous n'avons pu le construire en 
  l'inscrivant par la Géométrie dans une sphère donnée comme on l'a fait 
  pour les autres polyèdres.~»
Et pour cause, il s'agit là d'un problème solide (non résoluble à la règle
et au compas). Stévin n'a donc pas pu calculer le rapport entre l'arête et le
rayon de la sphère circonscrite. Quant à Descartes, il ne mentionne même pas 
ce polyèdre dans \cite{descartes}. Notons comme d'habitude $d$ l'arête du
polyèdre et $\rho$ le rayon de la sphère circonscrite. Rappelons que, 
si l'on pose
$$t=\frac{1}{3}\left(1+\sqrt[3]{19-3\sqrt{33}}
+\sqrt[3]{19+3\sqrt{33}}\right),$$
on a\footnote{Voir par exemple \cite{weissbach2002}. Il est d'ailleurs
remarquable qu'en 2002 les deux auteurs de cet article semblent n'avoir pu 
nulle part retrouver ce résultat \emph{exact} dans la littérature moderne sur 
les polyèdres. Ils écrivent~: «~In the literature, one can find relatively exact 
values for [$2\rho/d$] (or values coinciding with it up to a factor), or more 
rough estimates of [$2\rho/d$] are given. But one cannot find explanations 
how these estimates are obtained. This shows that a mathematical problem,
clearly formulated in Euclid's Book~XIII, was no longer in the center of
interests.~»}
$$2\rho=d\sqrt{\frac{3-t}{2-t}},$$
mais $\mathbb{Q}(t)$ n'est pas décomposable en un tour d'extensions
quadratiques sur $\mathbb{Q}$, c'est-à-dire que le 
rapport $\dfrac{2\rho}{d}$ n'est pas constructible à la règle et au compas.

Dans l'énoncé, Roth demande que $d=x_0+x_1$ où $x_0$ et
$x_1$ sont deux racines positives de l'équation
$$O_{2536}(x)=x^6+11x^5-\left(81+\frac{1}{2}\right)x^4-573x^3
+\left(2911+\frac{5}{16}\right)x^2-\left(4235+\frac{15}{16}\right)x
+\left(9756+\frac{1}{2}\right).$$
Cette équation est de la forme $P(x)=0$ avec
$$P(x)=x^6+11x^5-\left(81+\frac{1}{2}\right)x^4-573x^3
+\left(1644+\frac{5}{16}\right)x^2-\left(2969+\frac{15}{16}\right)x
+\left(9756+\frac{1}{2}\right).$$
On pose $x=\dfrac{z}{2}$~; alors
$$64P\left(\frac{z}{2}\right)=z^6+22z^5-326z^4-4584x^3+26309z^2-95038z
+624416.$$
On a $624416=2^5.13.19.79$, et on trouve un facteur quadratique de la forme
$z^2-yz+79$~:
$$64P\left(\frac{z}{2}\right)=(z+16)(z+26)(z^2+19)(z^2-20z+79).$$
Ainsi, le polynôme $P$ a seulement deux racines réelles positives
$\dfrac{10\pm\sqrt{21}}{2}$, et leur somme fait $d=10$~; on a donc
$$2\rho=10\sqrt{\frac{9-3t}{6-3t}}.$$

La solution donnée par Roth est-elle juste~? Son expression algébrique
avec les quatre racines cubiques ne laisse aucun doute quant au fait que
l'auteur de cette solution l'a trouvée en appliquant les formules de Cardan
à la construction de la sphère circonscrite. Normalisons numérateur, 
dénominateur et radicandes dans notre solution, à la manière de Roth, 
nous trouvons~:
$$10\sqrt{\frac{9-3t}{6-3t}}=\scriptstyle\sqrt{\frac{2666+\frac{2}{3}
-\sqrt[3]{703703703+\frac{19}{27}
-\sqrt{407407407407407407+\frac{11}{27}}}-\sqrt[3]{703703703+\frac{19}{27}
+\sqrt{407407407407407407+\frac{11}{27}}}}{16+\frac{2}{3}-\sqrt[3]{703
+\frac{19}{27}-\sqrt{407407+\frac{11}{27}}}-\sqrt[3]{703+\frac{19}{27}
+\sqrt{407407+\frac{11}{27}}}}}\displaystyle\simeq 26,8743.$$
C'est légèrement différent du résultat donné par Roth (arrondi au millième,
le résultat donné par Roth vaut $26,8728$). Tout se passe comme si l'auteur
de la solution donnée dans les \textit{Arithmetica philosophica} n'avait
gardé que les trois premiers chiffres des valeurs dans les radicandes,
et rempli un peu au hasard les autres rangs décimaux. Ainsi au dénominateur,
au lieu de
$$703+\frac{19}{27}=703,703703703\overline{703}$$
il aurait écrit
$$703+\frac{35}{54}=703,648148148\overline{148},$$
et au lieu de
$$407407+\frac{11}{27}=407407,407407\overline{407}$$
il aurait écrit
$$407329+\frac{215}{972}=407329,22\overline{119341563786008230452674897}.$$
Les nombres dans les radicandes du numérateur ne sont que des multiples
de ces deux nombres par des puissances de dix. On peut donc penser que
l'erreur commise n'est peut-être pas vraiment une faute, mais seulement un 
choix arbitraire quant à la précision retenue. Le calculateur aurait
négligé les décimales suivantes, ou bien les aurait calculées avec beaucoup
moins d'attention, sans se relire.

Il faut conclure que l'auteur de la solution connaissait en 1608 l'expression
exacte du rapport entre arête et diamètre de la sphère circonscrite au
\textit{cubus simus}.

\paragraph{III f. 189}
«~Soit un \textit{corpus irregulare} à 8 faces triangulaires et 18 faces carrées
équilatérales régulières, ainsi que 24 faces quadrilatères irrégulières (non
équilatérales), tel que chaque arête d'une face carrée s'attache à une face
quadrilatère irrégulière. Quand on le pose sur n'importe laquelle de ses
faces triangulaires, la face supérieure lui est toujours parallèle, et c'est
aussi un triangle. Chaque arête des faces carrées de ce \textit{corpus}
est en \textit{proportione dupla} avec chaque arête des faces triangulaires
\footnote{Probablement que $A$ est en \textit{proportione dupla} avec $B$ 
si $A=B\sqrt{2}$.}. Parmi les arêtes de chaque face quadrilatère irrégulière,
il y en a trois qui sont chacune une arête d'une face carrée, et la
quatrième est une arête d'une face triangulaire. L'arête des faces carrées de
ce \textit{corpus} est égale au produit des deux vraies racines
\textit{enneacositessaracontaoctogonales} de $1\cossvi+22\cossv-13\cossiv
-106\cossiii-1110\cossii-5657\cossi$. Combien le diamètre de la sphère
circonscrite mesure-t-il~? \emph{Solution}~: le diamètre de la sphère
mesure $\sqrt{288+\sqrt{41472}}$.~»

\paragraph{Commentaire}
Ce polyèdre irrégulier est obtenu par troncature à partir du
\textit{bistruncatum cubum primum} (objet de la question~XXI 
p.~\pageref{questioXXI} ci-dessus), les
arêtes étant coupées à la moitié. L'énoncé demande que l'arête des faces
carrées mesure $d=x_0\times x_1$ où $x_0$ et $x_1$ sont deux racines
réelles positives de l'équation
$$O_{948}(x)=x^6+22x^5-13x^4-106x^3-1110x^2-5657x.$$
Cette équation est de la forme $P(x)=0$ avec
$$P(x)=x^6+22x^5-13x^4-106x^3-1583x^2-5185x.$$
Mais la règle des signes, dite de Descartes\footnote{La règle des signes
était aussi connue des \textit{Rechenmeister}~; Faulhaber la publie dès 1622.},
montre aisément que ce polynôme n'a qu'une seule racine réelle positive (sans 
compter la racine nulle). Il doit donc y avoir une erreur d'énoncé.
\emph{Remarque}~: dans le polyèdre irrégulier décrit par Roth, on peut
montrer que $2\rho=2d\sqrt{2+\sqrt{2}}$. Or la solution donnée par Roth vaut
$2\rho=\sqrt{288+\sqrt{41472}}=12\sqrt{2+\sqrt{2}}$. Il faut donc que 
$d=x_0\times x_1=6$~: \textbf{pour corriger l'énoncé, on pourrait chercher
une équation qui ressemble à celle de l'énoncé, mais qui soit factorisable 
sur $\mathbb{Z}$ comme les précédentes, qui ait deux racines positives, et 
dont le terme constant soit divisible par 6.}

\paragraph{IV f. 189}
«~Soit un autre beau \textit{corpus irregulare} à 20 faces triangulaires et
12 faces pentagonales équilatérales régulières, ainsi que 30 faces
rectangulaires oblongues. Le rapport entre l'arête d'une face pentagonale
et l'arête d'une face triangulaire est comme le rapport de $\sqrt{72}$ à
$\sqrt{27-\sqrt{45}}$. Les arêtes des faces triangulaires sont égales aux
plus petites arêtes des faces rectangulaires, et les plus grandes arêtes
des faces rectangulaires sont égales aux arêtes des faces pentagonales.
L'arête d'une face pentagonale fait la somme des deux vraies racines
\textit{hexachilidiacosiaicosigonales} de $1\cossvi+27\cossv+3572\cossii
+22950-113\cossiv-669\cossiii-7353\cossi$. Combien le diamètre de la
sphère circonscrite à ce \textit{corpus} à soixante sommets mesure-t-il~?
\emph{Solution}~: $\sqrt{567+\sqrt{131220}}$.~»

\paragraph{Commentaire}
Le polyèdre ici en question n'est pas semi-régulier, mais il est très 
apparenté au rhombicosidodécaèdre de la question I p.~\pageref{questioI} 
ci-dessus. En fait il existe une famille continue de polyèdres convexes
$20[3]\ 12[5]\ 30[4]$ inscriptibles dans des sphères et obtenues par double
troncature à partir du dodécaèdre ou de l'icosaèdre, les faces triangulaires
et les faces pentagonales étant toutes régulières, et le rapport entre la 
largeur et la longueur des faces rectangulaires pouvant être quelconque.
Notons $d_1$ l'arête des faces pentagonales et $d_2$ l'arête des faces
triangulaires. Ici Roth pose 
$$\dfrac{d_1}{d_2}=\sqrt{\dfrac{72}{27-\sqrt{45}}}.$$
L'énoncé demande que $d_1=x_0+x_1$ où $x_0$ et $x_1$ sont deux racines
positives de l'équation
$$O_{6220}(x)=x^6+27x^5+3572x^2+22950-113x^4-669x^3-7353x.$$
Cette équation est de la forme $P(x)=0$ avec
$$P(x)=x^6+27x^5-113x^4-669x^3+463x^2-4245x+22950.$$
On a $22950=2.3^3.5^2.17$, et on trouve un facteur quadratique de la forme
$x^2-yx+17$~:
$$P(x)=(x+5)(x+30)(x^2+x+9)(x^2-9x+17).$$
Les deux seules racines réelles positives sont $\dfrac{9\pm\sqrt{13}}{2}$,
donc $d_1=9$. Si la valeur donnée par Roth comme solution s'avère exacte,
il faut donc que
$$\frac{(2\rho)^2}{d_1^2}=\frac{567+\sqrt{131220}}{81}=7+2\sqrt{5}\simeq 11{,}5.$$
Hélas, nous avons calculé que\footnote{\textbf{Reste à rédiger ce calcul, ou à trouver une référence qui contienne le résultat…}}
\[\frac{(2\rho)^2}{d_1^2}=3+3\varphi+(2+4\varphi)\frac{d_2}{d_1}+(2+\varphi)\left(\frac{d_2}{d_1}\right)^2\simeq 13{,}4.\]

\paragraph{V f. 190}
«~Soit un autre \textit{corpus irregulare}, beau, ingénieux et bien 
proportionné, à 12 faces pentagonales et 20 faces hexagonales équilatérales
régulières, ainsi que 60 faces triangulaires scalènes, tel qu'à chaque arête
des faces pentagonales est attaché une face triangulaire, et tel que, quand
on pose ce \textit{corpus} sur n'importe quelle face triangulaire, la
face supérieure lui est toujours parallèle, et c'est aussi un triangle.
Le rapport entre l'arête d'une face pentagonale et l'arête d'une face
hexagonale est comme le rapport de $\sqrt{18+\sqrt{180}}$ à $6$. L'arête
des faces pentagonales est égale à la plus petite arête des faces
triangulaires, mais les deux autres arêtes des faces triangulaires, plus
longues, sont chacune égale à l'arête des faces hexagonales. L'arête
des faces hexagonales mesure une unité de moins que le produit des deux
vraies racines \textit{myripentachilihexacosihebdomicontadygonales} de
$1\cossvi+8\cossv+7899\cossii+8415-56\cossiv-232\cossiii-8906\cossi$.
Combien le diamètre de la sphère circonscrite à ce \textit{corpus} à
quatre-vingt-dix sommets mesure-t-il~? \emph{Solution}~:
$\sqrt{5832+\sqrt{18895680}}$.~»

\paragraph{Commentaire}
Le polyèdre ici décrit n'est pas semi-régulier~; on l'obtient par troncature
du \textit{truncatum icosaedrum per laterum tertias} rencontré dans la
question~XXIV p.~\pageref{questioXXIV} ci-dessus, les arêtes étant coupées
à la moitié. Le polyèdre $12[5]\ 20[6]\ 60[3]$ ainsi obtenu est bien 
inscriptible dans une sphère mais ses faces triangulaires ne sont pas
équilatérales. Notons $d_1$ l'arête des faces pentagonales et $d_2$ l'arête
des faces hexagonales. Roth affirme que
$$\frac{d_1}{d_2}=\frac{\sqrt{18+6\sqrt{5}}}{6}.$$
Nous avons vérifié ce rapport\footnote{Voir note~\ref{rapport} ci-dessous.}. L'énoncé demande ensuite que $d_2=x_0x_1-1$ où $x_0$ 
et $x_1$ sont deux racines réelles positives de l'équation
$$O_{15672}(x)=x^6+8x^5-56x^4-232x^3+7899x^2-8906x+8415,$$
mais la solution donnée est erronée. Elle serait juste si l'énoncé demandait que \(d_2\) soit égal à \(x_0x_1+1\). C'est ce qu'on va montrer.
L'équation est de la forme $P(x)=0$ avec
$$P(x)=x^6+8x^5-56x^4-232x^3+64x^2-1072x+8415.$$
On a $8415=3^2.5.11.17$ et on trouve un facteur quadratique de la forme
$x^2-yx+17$ (déjà rencontré dans les questions précédentes)~:
$$P(x)=(x+5)(x+11)(x^2+x+9)(x^2-9x+17).$$
Les deux seules racines réelles positives sont $\dfrac{9\pm\sqrt{13}}{2}$, et
on a donc $d_2=x_0x_1+1=17+1=18$. Or la sphère circonscrite à ce polyèdre n'est autre
que la sphère tangente aux arêtes de l'icosaèdre d'arête \(3d\), où \(d_2=\dfrac{d\sqrt{3}}{2}\)~; mais\footnote{Voir \cite{coxeter1948} tabl.~I p.~293.} \(\dfrac{2\rho}{3d}=\varphi\). On en déduit \(\left(\dfrac{2\rho}{d_2}\right)^2=6(3+\sqrt{5})\) et
\[(2\rho)^2=18^2\times 6(3+\sqrt{5})=5832+\sqrt{18895680},\]
comme annoncé\footnote{\label{rapport}En passant, on remarque que, puisque \(d_1=\dfrac{\varphi d}{2}\), on a bien \(\dfrac{d_1}{d_2}=\dfrac{\varphi}{\sqrt{3}}=\dfrac{\sqrt{18+6\sqrt{5}}}{6}\), comme annoncé par Roth dans l'énoncé.}.

\end{document}